\documentclass[12pt]{article}

\def\ra{\rightarrow}

\def\ss{\subseteq}

\def\Im{\hbox{\rm Im}\,}

\def\Aut{{\rm Aut}}

 \def\HollowBox #1#2{{\dimen0=#1 \advance\dimen0 by -#2       
       \dimen1=#1 \advance\dimen1 by #2                       
        \vrule height #1 depth #2 width #2                    
        \vrule height 0pt depth #2 width #1                   
        \llap{\vrule height #1 depth -\dimen0 width \dimen1}%
       \hskip -#2                                             
       \vrule height #1 depth #2 width #2}}                   
 \def\BoxOpTwo{\mathord{\HollowBox{6pt}{.4pt}}\;}             

\def\endpf{\hfill $\BoxOpTwo$}

\font\teneufm=eufm10
\font\seveneufm=eufm7
\font\fiveeufm=eufm5
\newfam\eufmfam
\textfont\eufmfam=\teneufm
\scriptfont\eufmfam=\seveneufm
\scriptscriptfont\eufmfam=\fiveeufm

\newfam\msbfam
\font\tenmsb=msbm10 scaled \magstep1 \textfont\msbfam=\tenmsb
\font\sevenmsb=msbm7 scaled \magstep1   \scriptfont\msbfam=\sevenmsb
\font\fivemsb=msbm5 scaled \magstep1   \scriptscriptfont\msbfam=\fivemsb
\def\Bbb{\fam\msbfam \tenmsb}

\def\CC{{\Bbb C}}

\def\ZZ{{\Bbb Z}}

\newtheorem{theorem}{Theorem}[section]

\newtheorem{proposition}[theorem]{Proposition}

\newtheorem{example}[theorem]{EXAMPLE}

\makeindex

\usepackage{graphicx}

\usepackage{amsmath}

\begin{document}

\begin{center}
\huge \bf
Topological/Geometric Properties of the Orbit Accumulation Set\footnote{{\bf MR Classification
Numbers:}  32M05,32M99.}
\end{center}
\vspace*{.12in}

\begin{center}
\large Steven G. Krantz\footnote{Author supported in part
by the National Science Foundation and by the Dean of the Graduate
School at Washington University.}\end{center}
\vspace*{.15in}

\begin{center}
\today
\end{center}
\vspace*{.2in}

\begin{quotation}
{\bf Abstract:} \sl
We study the set of boundary orbit accumulation points of the automorphism
group action on a bounded domain in $\CC^n$.   Topological and geometric properties
of this set are derived.
\end{quotation}
\vspace*{.25in}

\setcounter{section}{-1}

\section{Introduction}

A {\it domain} $\Omega$ in $\CC^n$ is a connected, open set.  An {\it automorphism}
of $\Omega$ is a biholomorphic self-map.  The collection of automorphisms forms
a group under the binary operation of composition of mappings.   The topology
on this group is uniform convergence on compact sets, or the compact-open topology.
We denote the automorphism group by $\Aut(\Omega)$.

Although domains with {\it transitive automorphism group} are of some interest,
they are relatively rare (see [HEL]).   A geometrically more natural condition
to consider, and one that gives rise to a more robust and broader class of domains,
is that of having {\it non-compact automorphism group}.   Clearly a domain has
non-compact automorphism group if there are automorphisms $\{\varphi_j\}$ which
have no subsequence that converges to an automorphism.   The following proposition
of Henri Cartan is of particular utility in the study of these domains:

\begin{proposition} \sl
Let $\Omega \ss \CC^n$ be a bounded domain.  Then $\Omega$ has non-compact automorphism
group if and only if there are a point $X \in \Omega$, a point $P \in \partial \Omega$, and
automorphisms $\varphi_j$ of $\Omega$ such that $\varphi_j(X) \ra P$ as $j \ra \infty$.
\end{proposition}

\noindent We refer the reader to [NAR] for discussion and proof of Cartan's result.

A point $P$ in $\partial \Omega$ is called a {\it boundary orbit accumulation point} if there is an $X \in \Omega$ and
automorphisms $\varphi_j$ of $\Omega$ such that $\lim_{j\ra \infty} \varphi_j(X) = P$.

The inspiration for the present paper comes from [ISK1].  In that paper, we considered
the set $S(\Omega)$ of all boundary orbit accumulation points of a smoothly bounded, finite
type domain $\Omega$ in $\CC^n$ (see [KRA1] for the concept of finite type).   In particular, we showed that, if the set $S$ contains
at least three points, then it is perfect.   Also, if $S$ is not connected, then
it has uncountably many connected components.  Finally, the set $S$ is generically
contained in the set where the Levi form has minimal rank and where the type and the 
multiplicity functions are maximal.  

Our goal in this paper is to extend and generalize these results.  We also provide
examples to show that some of these statements are sharp.

\section{Some Examples}

We begin with some examples.   These will all be smooth deformations of the unit
ball in $\CC^2$.   

\begin{example} \rm
There is a bounded domain $\Omega$ in $\CC^2$ with Lipschitz boundary
with the property that $S(\Omega)$ has just two points.

To see this, let $B$ be the unit ball in $\CC^2$.  For $j \in \ZZ$, 
let 
$$
\varphi(z_1, z_2) = \left ( \frac{z_1 - 1/5}{1 - (1/5)z_1}, 
		 \frac{\sqrt{1 - (1/5)^2} z_2}{1 - (1/5)z_1} \right ) 
$$
and let $\varphi^j$ denote the  composition of $\varphi$ with
itself $j$ times (and of course $\varphi^j= (\varphi^{-1})^{|j|}$ when $j < 0$ and
$\varphi^0 = \hbox{id})$.

Then each $\varphi^j$ is an automorphism of $B$.  Furthermore,
$$
\lim_{j \ra +\infty} \varphi^j(z_1, z_2) = (-1, 0) \, ,
$$
uniformly on compact subsets of $(z_1, z_2)$ in $B$.   And
$$
\lim_{j \ra -\infty} \varphi^j(z_1, z_2) = (1, 0) \, ,
$$
uniformly on compact subsets of $(z_1, z_2)$ in $B$.

Now define 
$$
\alpha(z_1, z_2) = \left \{ \begin{array}{lcr}
                       (1/100 - |z_1|^2)^2 \cdot (1/100 - |z_2|^2)^2 & \hbox{if} & |(z_1, z_2)| \leq 1/10 \\
		              0               & \hbox{if} & |(z_1, z_2)| > 1/10 \, .
			    \end{array}
		   \right.
$$
We set
$$
\Omega' = \{(z_1, z_2) \in \CC^2: (|z_1|^2 + |z_2|^2 - 1) + \alpha(z_1 - i, z_2) < 0\} \, .
$$
Finally, let
$$
\Omega = \bigcap_{j = -\infty}^\infty \varphi^j(\Omega') \, .
$$

We see that $\Omega$ is the unit ball with infinitely many smooth dents which accumulate
at $(1,0)$ and $(-1,0)$.   Moreover, the automorphism group of $\Omega$ is just the
cyclic group $\{\varphi^j\}$ (see [LER], [GRK2] for the details).  
\end{example}

\begin{example} \rm
There is a bounded domain $\Omega$ in $\CC^2$ with Lipschitz boundary
with the property that $S(\Omega)$ has just one point.

To see this, let $B$ be the unit ball in $\CC^2$.  For $j \in \ZZ$, 
let 
$$
\varphi^j(z_1, z_2) = \left (  \frac{2i - j}{2i + j} \cdot \frac{z_1 + j/(j -2i)}{1 + z_1 j/(j + 2i)} , \frac{2i z_2/(j + 2i)}{1 + z_1 j/(j + 2i)} \right )  \, .
$$
The collection $\{\varphi^j\}_{j = -\infty}^\infty$ forms a cyclic subgroup of the automorphism group
of $B$.  One may see this by mapping $B$ to the Siegel upper half plane
${\cal U} = \{(w_1, w_2):  \Im w_1 > |w_2|^2\}$ and transferring the automorphims
to ${\cal U}$ (see [KRA2] for the details of this process).  
In fact, on ${\cal U}$, the automorphisms $\varphi^j$ become real translations
in the $w_1$ variable.

Then 
$$
\lim_{j \ra \pm \infty} \varphi^j(z_1, z_2) = (-1, 0) \, ,
$$
uniformly on compact subsets of $(z_1, z_2)$ in $B$.   
\end{example}

\section{Examples Based on the Bidisc}

The examples in this section are certainly of some utility, but they
are ancillary to the main points we are trying to make in this
paper.  We include them for completeness and for interest's sake.

\begin{example} \rm
There is a domain $\Omega \ss \CC^n$, modeled on the bidisc $D^2$, with
the property that $S(D)$ has dimension 0.
\end{example}

\begin{example} \rm
There is a domain $\Omega \ss \CC^n$, modeled on the bidisc $D^2$, with
the property that $S(D)$ has dimension 1.
\end{example}

\begin{example} \rm
There is a domain $\Omega \ss \CC^n$, modeled on the bidisc $D^2$, with
the property that $S(D)$ has dimension 2.
\end{example}

\begin{example} \rm
There is a domain $\Omega \ss \CC^n$, modeled on the bidisc $D^2$, with
the property that $S(D)$ has dimension 3.
\end{example}

\noindent {\bf Discussion of Examples 2.1--2.4:}  Let $P_0 \in \partial D^2$ be
given by $P_0 = (i,i)$.  Let 
$$
\varphi(z_1, z_2) = \left ( \frac{z_1 + 1/5}{1 + (1/5)z_1}, \frac{z_2 + 1/5}{1 + (1/5)z_2} \right ) \, .
$$
As usual, let $\varphi^j$ be the composition of $\varphi$ with itself $j$ times.  Then
$$
\lim_{j \ra \pm \infty} \varphi^j(P_0) = (\pm 1, \pm 1) \, .
$$
If we set up a ``dent'' at $P_0$, just as in Examples 1.1 and 1.2, and propagate the dent
using the automorphisms $\varphi^j$, then we produce a domain with automorphism group
a cyclic group and two orbit accumulation points.  This is a zero-dimensional $S(\Omega)$.

Now let us modify the discussion in the preceding paragraph by defining the disc automorphism
$$
\lambda(\zeta) = \frac{\zeta + 1/5}{1 +(1/5)\zeta} 
$$
and then considering the automorphisms
$$
\psi_{j, \rho}(z_1, z_2) = \left ( \lambda^j(z_1), \rho^j(z_2) \right ) \, ,
$$
where $\lambda^j$ as usual denote the $j$-fold composition of $\lambda$ with itself
and $\rho$ denotes {\it any} nontrivial M\"{o}bius transformation $\zeta \mapsto (\zeta - a)/(1- \overline{a}\zeta)$ of the disc $D$.

As usual, we construct a small dent at the boundary point $(i,0)$ and we propagate
it using the automorphisms $\psi_{j, \rho}$ for all integers $j$ and all automorphisms
$\rho$ of the unit disc.   Then it is easy to see that the boundary orbit accumulation
points of the resulting domain $\Omega$ are points of the form
$$ 
\{1\} \times \partial D \qquad \hbox{and} \qquad \{-1\} \times \partial D \, .
$$
We see that, according to this construction, the dimension of $S(\Omega)$ is 1.

Now we modify the construction in the last two paragraphs by letting
$$
\mu_j(z_1, z_2) = \left ( \lambda^j(z_1), z_2 \right ) \, ,
$$
We construct a domain with dents as usual---a modification of the bidisc $D^2$---using these 
automorphisms.   It is easy to see that, for the resulting $\Omega$, 
the set $S(\Omega)$ is 
$$
\bigl ( \{-1\} \times D \bigr ) \ \cup \ \bigl ( \{1\} \times D \bigr ) \, .
$$
Thus $S(D)$ is of dimension 2.

Finally, if we consider just the usual bidisc $D^2$ with its full
automorphism group, then the set $S(D^2)$ is the entire topological
boundary $\partial D \times D \cup D \times \partial D$.  Thus
$S(D^2)$ has dimension 3.

\section{General Principles}

Now we wish to make some general assertions about the geometric nature of the
orbit accumulation set $S(\Omega)$.

\begin{proposition} \sl
Let $\Omega \ss \CC^n$ be a smoothly bounded domain of finite type
in the sense of Catlin/D'Angelo/Kohn (see [KRA1] for details).  
Then it is impossible for $S(\Omega)$ to contain a relatively
open subset $U$ of $\partial \Omega$ unless $\Omega$ is biholomorphic to
the unit ball $B$.
\end{proposition}
{\bf Proof:}   This assertion is particularly easy to see in complex
dimension 2.  For if a relatively open subset $U$ of the boundary has all points
of geometric type greater than 2, then $U$ is foliated by 1-dimensional
complex analytic varieties (see [FRE]).\footnote{One could at this point invoke the theorem of Kim [KIM]---at least
in the convex case---to see that $\Omega$ must be biholomorphic to the bidisc.  That would be a contradiction.}
So every point of $U$ is in fact of infinite
type and that is a contradiction.   Thus any relatively open subset $U$ of $\partial \Omega$
will contain points of type 2.  Of course these points must be pseudoconvex by the result of [GRK1].
Hence the points are strongly pseudoconvex.  But then the theorem of Bun Wong and Rosay
(see [WON] and [ROS]) tells us that the domain $\Omega$ is biholomorphic to the unit ball
and hence $S(\Omega)$ is the entire sphere.  That proves the result.

In $\CC^n$ with $n > 1$ we must instead use the structure theorem of Catlin [CAT]
which tells us that the set of finite type points with the type greater than 2 
forms a subset of the boundary of codimension at least 1.   Hence any relatively
open subset of the boundary will contain pseudoconvex points of type 2, and the
argument is finished as before.
\endpf
\smallskip \\

\begin{proposition} \sl
Let $\Omega \ss \CC^n$ be a smoothly bounded domain of finite type
in the sense of Catlin/D'Angelo/Kohn (see [KRA1] for details).  
Then it is impossible for $S(\Omega)$ to contain a subset $E$ of $\partial \Omega$ having
positive $(2n-1)$-dimensional Hausdorff measure unless $\Omega$ is biholomorphic to
the unit ball $B$.
\end{proposition}
{\bf Proof:}  Reasoning as in the second part of the proof of the preceding proposition,
we see that the set $E$ must contain a strongly pseudoconvex point.  It follows then
from the Bun Wong/Rosay theorem that the domain must be biholomorphic to the unit
ball.
\endpf
\smallskip \\

\section{Closing Thoughts}

There is no Riemann mapping theorem in several complex variables.
So we seek other means to compare and contrast domains up to 
biholomorphic equivalence.  Certainly the automorphism group
is one invariant that has proved to be both flexible and useful.

One avenue that has been explored to a great extent is the Levi
geometry of boundary orbit accumulation points (see [GRK3] and [ISK2]).
The nature of the boundary orbit accumulation set is a much newer
avenue of inquiry, and there remain many questions to be answered
in this venue.  We hope to explore some of them in future papers.

\newpage

\null \vspace*{-1in}

\noindent {\Large \sc References}
\bigskip  \\

\begin{enumerate}

\item[{\bf [BEB]}]  S. Bell and H. Boas, Regularity of the
Bergman projection in weakly pseudoconvex domains, {\it Math.\
Ann.} 257(1981), 23--30.

\item[{\bf [CAT]}]  D. Catlin, Boundary invariants of
pseudoconvex domains, {\it Ann.\ of Math.} 120(1984), 529--586

\item[{\bf [CHM]}]  S.-S. Chern and J. Moser, Real hypersurfaces in complex
manifolds, {\em Acta Math.} 133(1974), 219-271.

\item[{\bf [FRE]}]  M. Freeman, Local complex foliation of real
submanifolds, {\it Math. Ann.} 209(1974), 1-30.

\item[{\bf [GRK1]}]  R. E. Greene and S. G. Krantz, Techniques
for studying automorphisms of weakly pseudoconvex domains,
{\it Several Complex Variables} (Stockholm, 1987/1988),
Math.\ Notes 38, Princeton University Press, Princeton, NJ,
1993, 389--410.

\item[{\bf [GRK2]}]  R. E. Greene and S. G. Krantz, Invariants
of Bergman geometry and the automorphism groups of domains in
$\CC^n$, {\it Geometrical and Algebraical Aspects in Several
Complex Variables} (Cetraro, 1989), 107--136, Sem.\ Conf., 8,
EditEl, Rende, 1991.

\item[{\bf [GRK3]}] R. E. Greene and S. G. Krantz,
Biholomorphic self-maps of domains, {\it Complex Analysis II}
(C. Berenstein, ed.), Springer Lecture Notes, vol. 1276, 1987,
136-207.

\item[{\bf [HEL]}]  S. Helgason, {\it Differential Geometry and
Symmetric Spaces}, Academic Press, New York, 1962.

\item[{\bf  [ISK1]}]  A. Isaev and S. G. Krantz, On the boundary
orbit accumulation set for a domain with non-compact
automorphism group, {\it Mich.\ Math.\ Jour.} 43(1996),
611-617.

\item[{\bf [ISK2]}] A. Isaev and S. G. Krantz, Domains with
non-compact automorphism group: a survey, {\it Adv.\ Math.}
146(1999), 1--38.

\item[{\bf [KIM]}]  K.-T. Kim, Domains in $\CC^n$ with a
piecewise Levi flat boundary which possess a noncompact
automorphism group, {\it Math. Ann.} 292(1992), 575--586.

\item[{\bf [KIK]}] K.-T. Kim and S. G. Krantz, Complex scaling
and domains with non-compact automorphism group, {\it Illinois
Journal of Math.} 45(2001), 1273--1299.
			  
\item[{\bf [KRA1]}] S. G. Krantz, {\it Function Theory of
Several Complex Variables}, $2^{\rm nd}$ ed., American
Mathematical Society, Providence, RI, 2001.

\item[{\bf [NAR]}] R. Narasimhan, {\it Several Complex
Variables}, University of Chicago Press, Chicago, 1971.

\item[{\bf [KRA2]}] S. G. Krantz, {\it Explorations in Harmonic
Analysis, with Applications in Complex Function Theory and the
Heisenberg Group}, Birkh\"{a}user Publishing, Boston, 2009. \\

\item[{\bf [LER]}] L. Lempert and L. A. Rubel, An independence
result in several complex variables, {\it Proc.\ Amer.\ Math.\
Soc.} 113(1991), 1055--1065.

\item[{\bf [NAR]}] R. Narasimhan, {\it Several Complex
Variables}, University of Chicago Press, Chicago, 1971.

\item[{\bf [ROS]}] J.-P. Rosay, Sur une characterization de la
boule parmi les domains de $\CC^n$ par son groupe
d'automorphismes, {\it Ann. Inst. Four. Grenoble} XXIX(1979),
91-97.

\item[{\bf [WON]}] B. Wong, Characterizations of the ball in
$\CC^n$ by its automorphism group, {\it Invent. Math.}
41(1977), 253-257.

\end{enumerate}
\vspace*{.95in}

\small

\noindent \begin{quote}
Department of Mathematics \\
Washington University in St. Louis \\
St.\ Louis, Missouri 63130 \\ 			   
{\tt sk@math.wustl.edu}
\end{quote}

\end{document}